\documentclass[titlepage,11pt]{article}
\usepackage{latexsym}
\usepackage{amsfonts}
\usepackage{amsmath}
\newcommand\blackslug{\hbox{\hskip 1pt \vrule width 4pt height 8pt depth 1.5pt
        \hskip 1pt}}
\newcommand\bbox{\hfill \quad \blackslug \bigbreak}
\def\d{\hbox{-}}
\def\c{\hbox{-}\cdots\hbox{-}}
\def\l{,\ldots,}

\title{Induced subgraphs of graphs with large chromatic number.\\
XI. Orientations}
\author{Maria Chudnovsky\thanks{Supported by NSF grant DMS-1550991 and
US Army
Research Office Grant W911NF-16-1-0404.}\\
Princeton University, Princeton, NJ 08544
\\
\\
Alex Scott\thanks{Supported by a Leverhulme Trust Research Fellowship.}\\
Oxford University, Oxford, UK
\\
\\
Paul Seymour\thanks{Supported by ONR grant N00014-14-1-0084 and 
NSF grant DMS-1265563.}\\
Princeton University, Princeton, NJ 08544}

\date{February 6, 2017; revised \today}

\newtheorem{thm}{}[section]

\newcommand{\Proof}{\noindent{\bf Proof.}\ \ }

\begin{document}
\maketitle
\begin{abstract}
Fix an oriented graph $H$, and let $G$ be a graph with bounded clique number and very large chromatic number. 
If we somehow orient its edges, must there be an induced subdigraph isomorphic to $H$?
Kierstead and R\"odl~\cite{kiersteadrodl} raised this question for two specific kinds of digraph $H$:
the three-edge path, with the first and last edges 
both directed towards the interior; and stars (with many edges directed out and many directed in).
Aboulker et al.~\cite{aboulker} subsequently conjectured that the answer is affirmative in both cases.
We give affirmative answers to both questions.

\end{abstract}

\section{Introduction}

All graphs in this paper are finite and simple.
If $G$
is a graph, $\chi(G)$ denotes its chromatic number, and $\omega(G)$ denotes its clique number, that is, the cardinality
of the largest clique of $G$. This paper is concerned with the digraphs that can be obtained by orienting the edges of a graph,
and in particular, digraphs in this paper have no ``antiparallel'' pairs of edges, that is, no directed cycles of 
length two, as well as no loops or parallel edges.  If $G$ is a digraph, $G^*$ means the underlying graph.
For a digraph $G$, we say $u$ is $G$-adjacent {\em to $v$} or {\em from $v$} 
to indicate the direction of the edge between $u$ and $v$, and $u$ is $G^*$-adjacent {\em with $v$} to mean adjacency in $G^*$.
The {\em chromatic number} $\chi(G)$ and {\em clique number} $\omega(G)$ of a digraph $G$ mean the corresponding quantities for $G^*$.

Let $H$ be a graph.  We say that $H$ is {\em $\chi$-bounding} if there is a function $f$ such that $\chi(G)\le f(\omega(G))$
for every graph $G$ not containing $H$ as an induced subdigraph.  If $H$ is $\chi$-bounding then it cannot contain a cycle, since 
(as shown by Erd\H os \cite{erdos}) there are graphs with large girth and large chromatic number.  Thus the only possible $\chi$-bounding graphs are forests.
The Gy\'arf\'as-Sumner conjecture~\cite{gyarfastree, sumner} asserts that every forest is $\chi$-bounding.  Despite considerable work, the conjecture is only known for
comparatively few families (see \cite{gyarfas,gst,kp,kz,scott,newbrooms,caterpillars}).

Now let $H$ be an oriented graph. When is there a function $f$ such that $\chi(G)\le f(\omega(G))$
for every oriented graph $G$ not containing $H$ as an induced subdigraph?  As in the graph case, we call a digraph $H$ with this property {\em $\chi$-bounding}.
Then every $\chi$-bounding digraph $H$ is an oriented forest, because we can take
$G$ to be any orientation of a graph with large girth and large chromatic number. 
However, for digraphs it is {\em not} the case that $H$ is $\chi$-bounding whenever $H^*$ is a forest.
Indeed, this is false even for digraphs $H$ such that $H^*$ is a three-edge path.
 
There are four ways to orient the edges of a three-edge path, up to reversing the path,
and we denote the corresponding digraphs by 
$$\rightarrow  \rightarrow \rightarrow, \rightarrow \leftarrow \rightarrow, \rightarrow \leftarrow \leftarrow, \leftarrow \rightarrow \rightarrow$$
with the natural meaning. Then
\begin{itemize}
\item Kierstead and Trotter~\cite{kiersteadtrotter} showed that $\rightarrow  \rightarrow \rightarrow $ is not $\chi$-bounding, by constructing triangle-free graphs
with arbitrarily large chromatic number together with a suitable orientation;
\item Gy\'arf\'as \cite{gyarfasproblem} noted that  
$\rightarrow \leftarrow \rightarrow$ is not $\chi$-bounding: let $D$ be the natural orientation of the shift graph on pairs, so $D$ has vertex set $[n]^{(2)}$, with edges
from $\{i,j\}$ to $\{j,k\}$ whenever $i<j<k$.  Then $D$ is triangle-free, has large chromatic number, and does not contain an induced copy of $\rightarrow \leftarrow \rightarrow$.
\end{itemize}
The two remaining orientations are equivalent under reversing all edges, so both or neither are $\chi$-bounding.  Thus it is enough to consider $\rightarrow \leftarrow \leftarrow$.

Oriented graphs with no induced $\rightarrow \leftarrow \leftarrow$ have been considered by several authors.
In the special case of {\em acyclic} orientations, Chv\'at\'al \cite{chvatal} showed that if $G$ is an acyclic oriented graph with no induced copy of 
$\rightarrow \leftarrow \leftarrow$ then $G$  is perfect, and so $\chi(G)=\omega(G)$.  Kierstead and R\"odl \cite{kiersteadrodl} asked whether  $\rightarrow \leftarrow \leftarrow$ is $\chi$-bounding,
and showed that the class of oriented graphs with no induced copy of $\rightarrow \leftarrow \leftarrow$ and no cyclic triangle is  $\chi$-bounded.
Aboulker et al.~\cite{aboulker} conjectured that $\rightarrow \leftarrow \leftarrow$ is in fact $\chi$-bounding, and proved some further special cases.
Our first main result resolves the question.
\begin{thm}\label{mainthm1}
The digraphs $\rightarrow \leftarrow \leftarrow$ and $\leftarrow \rightarrow \rightarrow$ are $\chi$-bounding.
\end{thm}
This is proved in the next section.

Which forests have the property that every orientation is $\chi$-bounding?  The results of Gy\'arf\'as and of Kierstead and Trotter mentioned above show that such a forest cannot contain a three-edge path, and so every component must be a star.

A digraph $H$ is an {\em oriented star} if $H^*$ is a star, that is, isomorphic to the complete bipartite graph $K_{1,t}$ for some $t\ge 0$.
As with paths, there have been several previous results on the chromatic number of graphs with a forbidden oriented star.
Gy\'arf\'as \cite{gyarfasproblem} asked whether, for every oriented star $H$, the class of acyclic oriented graphs with no induced $H$ is $\chi$-bounded.
Kierstead and R\"odl \cite{kiersteadrodl} proved the stronger result that the class of oriented graphs with no induced $H$ and no cyclic triangle is $\chi$-bounded; 
they further asked whether every oriented star is $\chi$-bounding.  
Aboulker et al.~\cite{aboulker} conjectured that oriented stars are indeed $\chi$-bounding, and showed that for every oriented star $H$ the
class of oriented graphs with no induced $H$ and no transitive triangle has bounded chromatic number (note that every orientation of $K_4$ has a transitive triangle, so if $G$ has no transitive triangle then $\omega(G)$ is at most 3).
Our second main result  answers this question.
\begin{thm}\label{stars}
Every oriented star is $\chi$-bounding.
\end{thm}
It is easy to prove this for stars in which every edge is directed away from the centre, or every 
edge is directed towards the centre, but the case when there are edges of both types is more difficult.
It follows from \ref{stars} that if $F$ is a forest such that every component is a star then every orientation of $F$ is $\chi$-bounding.

\section{An oriented three-edge path}

If $X\subseteq V(G)$, $G[X]$ denotes
 the subgraph or subdigraph induced on $X$, and we write $\chi(X)$ for $\chi(G[X])$ when there
is no danger of ambiguity.
If $G$ is a digraph and $v\in V(G)$, we denote the set of vertices with distance at most $r$ (in $G^*$) from $v$ by $N^{r}[v]$
or $N^r_G[v]$, 
and the set with distance
exactly $r$ by $N^{r}(v)$. We denote by $\chi^{r}(G)$ the maximum of $\chi(N^{r}[v])$ over all $v\in V(G)$ 
(or zero for the null digraph.) 

In this section we prove our first 
main result, that $\rightarrow \leftarrow \leftarrow$ is $\chi$-bounding. 
In fact, with a very little extra work we can prove a stronger statement, which
we now explain. A {\em hole} in a graph is an induced cycle of length at least four, and when $G$ is a digraph, by a ``hole'' of $G$ we
mean an induced subgraph $C$ such that $C^*$ is a hole of $G^*$. By a {\em long} hole we mean (just in this paper)
a hole of length at least five. A hole of a digraph $C$ is 
\begin{itemize} 
\item {\em directed} if each of its vertices has outdegree one in $C$; 
\item {\em alternating} if each of its vertices has outdegree two or zero in $C$ (and therefore $C$ has even length); and
\item {\em disoriented} if it is neither directed nor alternating.
\end{itemize}
It is easy to see that if some long hole of $G$ is disoriented, then $G$ contains 
$\rightarrow \leftarrow \leftarrow$ as an induced subdigraph.
(Some two consecutive edges of $G$ make a two-edge directed path, but $C$ is not a directed
cycle; grow the path to a maximal directed path of $C$ and look at its ends.)
Thus the following theorem implies that $\rightarrow \leftarrow \leftarrow$ is $\chi$-bounding.
(A useful feature of this strengthening is that now we are proving something invariant under reversing all edges of $G$,
which reduces the case analysis.)

\begin{thm}\label{holes}
For all $\kappa$ there exists $c$ such that if $G$ is a digraph with $\omega(G)\le \kappa$ and $\chi(G)>c$ then some long hole
of $G$ is disoriented.
\end{thm}
\Proof
We proceed by induction on $\kappa$; thus we may assume that $\chi(J)\le \tau$ for every digraph $J$ with 
with $\omega(J)<\kappa$ and no disoriented long hole. Let $c=2(3\tau)^5$; we claim that $c$ satisfies the theorem.
Let $G$ be a digraph with $\omega(G)\le \kappa$
and with no disoriented long hole. 
\\
\\
(1) {\em For each vertex $z$ and integer $r\ge 1$, $\chi(N^r(z))\le 3\tau\chi(N^{r-1}(z))$. Consequently
$\chi(N^r(z))\le \tau(3\tau)^{r-1}$ and $\chi(N^s(z))\le (3\tau)^{s-r}\chi(N^r(z))$ for all $s\ge r$ .}
\\
\\
Let us write $L_i$ for $N^i(z)\;(i\ge 0)$. Since $\omega(G[L_1])<\kappa$ the result holds if $r=1$, so we assume $r\ge 2$.
Let $I\subseteq L_{r-1}$ be stable. Let $I_1$ be the set of vertices in $I$
with no in-neighbours in $L_{r-2}$; $I_2$ the set with no out-neighbours in $L_{r-2}$;
and $I_3=I\setminus (I_1\cup I_2)$. Let $J_i$ be the set of vertices in $L_r$ with a neighbour in $I_i$
for $i = 1,2,3$. Let $i\in \{1,2,3\}$, and suppose that $\omega(G[J_i])=\kappa$.
Choose a clique $K$ of $G[J_i]$ with cardinality $\kappa$, and take a minimal subset $I_0$ of $I_i$ such that every vertex in $K$
has a neighbour in $I_0$. Since $\omega(G)=|K|$, it follows that $|I_0|\ge 2$; choose distinct $v_1,v_2\in I_0$. From the minimality
of $I_0$, there exists $u_1\in K$ $G^*$-adjacent with $v_1$ and not with $v_2$, and $u_2$ $G^*$-adjacent with $v_2$ and not with $v_1$.
Thus $v_1\d u_1\d u_2\d v_2$ is an induced path of $G^*$.
By reversing all edges of $G$ if necessary (this is legitimate since what we are proving is invariant under this reversal)
we may assume that $i\in \{1,3\}$. 

Since $G^*[L_0\cup\cdots\cup L_{r-2}]$ is connected, there is an induced path of $G^*$ joining $v_1,v_2$ with interior in this set,
and its union with  $v_1\d u_1\d u_2\d v_2$ is a hole. We may assume this hole is either directed or alternating in $G$, 
and in either case
exactly one of the edges $u_1v_1, u_2v_2$ of $G^*$ is oriented in $G$ from $L_r$ to $L_{r-1}$. Consequently 
we may assume that $v_1u_1$
and $u_2v_2$ are edges of $G$. Since $i\in \{1,3\}$, both $v_1,v_2$ have out-neighbours in $L_{r-2}$, say $w_1,w_2$ respectively.
If $w_1$ is $G^*$-adjacent with $v_2$, then adding $w_1$ to $u_1,v_1,v_2,u_2$ gives a hole of length five that is not directed,
a contradiction; so $w_1,v_2$ are $G^*$-nonadjacent. In particular $w_1\ne w_2$, so $r\ge 3$. If $v_1,w_2$ are $G^*$-nonadjacent, there is an 
induced path of $G^*$ between $w_1$ and $w_2$ with interior in $L_0\cup\cdots\cup L_{r-3}$, and its union with 
$w_1\d v_1\d u_1\d u_2\d v_2\d w_2$ yields a disoriented long hole of $G$, a contradiction; so $v_1,w_2$ are $G^*$-adjacent. This provides a hole
of length five, which is therefore directed; so $u_1u_2$ and $w_2v_1$ are edges of $G$. In particular $v_1$ has both an in-neighbour
and an out-neighbour in $L_{r-2}$, and so $i=3$, and therefore $v_2$ has an in-neighbour $x_2$ say in $L_{r-2}$. Since
the path $u_2v_2x_2$ is not directed, it follows that $x_2,v_1$ are $G^*$-nonadjacent, and in particular $x_2\ne w_1$. Join
$w_1,x_2$ by an induced path with interior in $L_0\cup\cdots\cup L_{r-3}$; then the union of this with the path
$w_1\d v_1\d u_1\d u_2\d v_2\d x_2$ yields a disoriented long hole, a contradiction. This proves that $\omega(G[J_i])<\kappa$,
and so $\chi(J_i)\le \tau$. Consequently $\chi(J_1\cup J_2\cup J_3)\le 3\tau$. Applying this to each colour class of
a $\chi(L_{r-1})$-colouring of $G[L_{r-1}]$, we deduce the 
first assertion of (1). The second
follows from the first by induction on $r$, since $\omega(G[L_1])<\kappa$ and so $\chi(L_1)\le \tau$; 
and the third follows from the first
by induction on $s-r$. This proves (1).

\bigskip

Suppose that $\chi(G)>c=2(3\tau)^5$. We may assume that $G^*$ is connected; choose a vertex $z$, and let $L_i=N^i(z)$ for all $i\ge 0$.
Choose $s$ such that $\chi(L_s)\ge \chi(G)/2$. Since $\chi(G)>2(3\tau)^5$, it follows that $\chi(L_s)>(3\tau)^{5}$ and so $s\ge 6$
by (1). Let $S$ be the vertex set of a component of $G[L_s]$ with maximum chromatic number. 
Let $r=s-4$, and choose $R\subseteq L_{r}$ minimal such that every vertex in $S$ is joined to a vertex in $R$ by a path 
in $G^*$
of length $4$. Let $G'=G\setminus (L_{r}\setminus R)$. Thus $N^{r}_{G'}(z) = R$, and $S\subseteq N^s_{G'}$. By the third assertion 
of (1)
applied to $G'$, $\chi(N^s_{G'})\le (3\tau)^4\chi(N^r_{G'})$, and since $\chi(N^s_{G'})>3(3\tau)^{4}$
it follows that $\chi(R)>2$ (indeed, $\chi(R)>3$). 

If $a\in R$ and $v\in S$, we say that $a$ is an {\em ancestor}
of $v$ if there is a path of
$G^*$ between $v$ and $a$ of length $4$.
From the minimality of $R$,
for each $a\in R$ there is a vertex $v$ in $S$ such that $a$ is its unique ancestor; let $P_a$ be a path
between $a$ and some such $v$ of length $4$.

Let $R_1$ be the set of vertices $a\in R$ such that the edge of $P_a$ incident with $a$ has head $a$, and $R_2$
the set for which this edge has tail $a$. Since $\chi(R)>2$, not both $R_1,R_2$ are stable, and by reversing all edges if necessary,
we may assume that $R_1$ is not stable. Let $a_r,b_r\in R_1$ be $G^*$-adjacent. Let the vertices of $P_{a_r}$ be $a_r\d a_{r+1}\c a_s$
in order, and let those of $P_{b_r}$ be $b_r\d b_{r+1}\c b_s$ in order. Since $b_r$ is the unique ancestor of $b_s$,
$b_i$ is $G^*$-nonadjacent with $a_{i-1}$ for $r+1\le i\le s$, and similarly $a_i$ is $G^*$-nonadjacent with $b_{i-1}$ for $r+1\le i\le s$.
In particular, $a_i\ne b_i$ for $r\le i\le s$. (However, $a_i, b_i$ may be adjacent in $G^*$.)
\\
\\
(2) {\em There is an induced path of $G^*$ between $a_{r+2}$ and $b_r$ whose interior contains no neighbours of $a_{r+1}, b_{r+1}, b_{r+2}$,
and an induced path between $b_{r+2}$ and $a_r$ whose interior contains no neighbours of $a_{r+1}, b_{r+1}, a_{r+2}$.}
\\
\\
Since $\chi(S)>3\chi^{5}(G)$, there is a vertex $v\in S$ with distance at least $6$ from each of $a_{r+1}, b_{r+1}, b_{r+2}$.
Let $u\in R$ be an ancestor of $v$, and let $P$ be a path of length $4$ between $u$ and $v$. Thus none of $a_{r+1}, b_{r+1}, b_{r+2}$
have neighbours in $V(P)$. Also there is an induced path between $u$ and $b_r$ with interior in $L_0\cup\cdots\cup  L_{r-1}$, an induced path
between $a_s$ and $v$ with interior in $S$, and the path $a_{r+2}\c a_s$. The union of these paths gives a path of $G^*$ (not necessarily 
induced) between $a_{r+2}$ and $b_r$, and so there is an induced path of $G^*$ using a subset of the same vertices between 
$a_{r+2}$ and $b_r$.
None of $a_{r+1}, b_{r+1}, b_{r+2}$ has a neighbour in any of these paths except $a_{r+2}, b_r$, and this proves the first statement.
The second follows by symmetry. This proves~(2).

\bigskip

Since $G^*[S]$ is connected, there is an induced path of $G^*$ between $a_{r+1}$ and $b_{r+1}$ with interior in 
$L_{r+2}\cup L_{r+3}\cup S$, and 
since the union of this path with the path $a_{r+1}\d a_r\d b_r\d b_{r+1}$ does not give a disoriented long hole, it follows
that $a_{r+1}, b_{r+1}$ are $G^*$-adjacent. By exchanging $a_r,b_r$ if necessary, we may assume that the edge 
$a_{r+1}b_{r+1}$ has head $b_{r+1}$. If the edge $a_{r+2}a_{r+1}$ has head $a_{r+2}$, the path $b_r\d b_{r+1}\d a_{r+1}\d a_{r+2}$
together with the first path of (2) gives a disoriented long hole, a contradiction. So $a_{r+2}a_{r+1}$ has head $a_{r+1}$. If
the edge $b_{r+2}b_{r+1}$ has head $b_{r+2}$, then the path $a_r\d a_{r+1}\d b_{r+1}\d b_{r+2}$ together with
the second path of (2) gives a disoriented long hole; so $b_{r+2}b_{r+1}$ has head $b_{r+1}$.
There is an induced path joining $a_{r+2}, b_{r+2}$ with interior in $L_{r+3}\cup S$,
and its union with $a_{r+2}\d a_{r+1}\d b_{r+1}\d b_{r+2}$ does not give a disoriented long hole; so $a_{r+2}, b_{r+2}$ are $G^*$-adjacent.
If the edge $a_{r+2}b_{r+2}$ has head $a_{r+2}$, the union of the path $a_{r+2}\d b_{r+2}\d b_{r+1}\d b_r$ with the first path of (2)
gives a disoriented long hole; while if $a_{r+2}b_{r+2}$ has head $b_{r+2}$, the union of $b_{r+2}\d a_{r+2}\d a_{r+1}\d a_r$
with the second path of (2) gives a disoriented long hole. This proves~\ref{holes}.~\bbox

\section{Oriented stars}

Now we turn to the proof of \ref{stars}. 
If $v$ is a vertex of a digraph $G$, $N^+(v)$ denotes the set of out-neighbours of $V$, and $N^-(v)$ denotes the set of in-neighbours.
Let us say a digraph $G$ is {\em $\lambda$-spread} if for every vertex $v$ of $G$,
and for all $A\subseteq N^+(v)$ and $B\subseteq N^-(v)$ with $|A|=|B|=\lambda$, some vertex of $A$ is $G^*$-adjacent with some vertex of $B$.
If $S$ is an oriented star and $\kappa\ge 0$, choose $\lambda$ such that every 
graph with at least $\lambda$ vertices has either a clique of cardinality
more than $\kappa$ or a stable set of cardinality $|V(S)|$. It follows then that every digraph $G$ with $\omega(G)\le \kappa$
not containing $S$ as an induced subdigraph is $\lambda$-spread. (For otherwise, with $v,A,B$ as above, 
there are stable subsets of $A, B$ each of 
cardinality $|V(S)|$, and an appropriate  subset of $A\cup B$
together with $\{v\}$ induces $S$, a contradiction.) It is more convenient to replace the hypothesis that 
$G$ does not contain $S$ as an induced subdigraph with the hypothesis that $G$ is $\lambda$-spread. Thus, now we need to prove:
for all $\kappa,\lambda\ge 0$, if $G$ is a $\lambda$-spread digraph with $\omega(G)\le \kappa$ then the 
chromatic number of $G$ is bounded
by some function of $\kappa,\lambda$. We will prove this by induction on $\kappa$, for fixed $\lambda$; 
so we will assume (throughout this
section) that $\kappa, \lambda$ and $\tau$ are fixed integers satisfying
\begin{itemize}
\item  $\kappa\ge 2$, $\lambda\ge 0$, and $\chi(J)\le \tau$ for 
every $\lambda$-spread digraph $J$ with $\omega(J)<\kappa$.
\end{itemize}
If $G$ is a graph and $A,B\subseteq V(G)$, we say $A$ is {\em $G$-complete with $B$} if $A\cap B=\emptyset$ and
every vertex in $A$ is $G$-adjacent with every vertex in $B$. If $G$ is a digraph, we say $A$ is {\em $G$-complete to $B$} and
{\em $B$ is complete from $A$}
if $A\cap B=\emptyset$ and
every vertex in $A$ is $G$-adjacent
to every vertex of $B$.
We need the following, which is an easy application of Ramsey's theorem~\cite{ramsey} and its bipartite version~\cite{beineke},
and we omit its proof:

\begin{thm}\label{bigramsey}
For all $k,m$ there exists $n\ge 0$ with the following property. Let $A_1\l A_n, B_1\l B_n$
be pairwise disjoint subsets of the vertex set of a graph $G$, each of cardinality $m$. 
Then either 
\begin{itemize}
\item there exist $A\subseteq A_1\cup \cdots\cup A_n$ and $B\subseteq B_1\cup \cdots\cup B_n$ with $|A|=|B|=\lambda$
such that no vertex in $A$ has a neighbour in $B$, or
\item there exist $I, J\subseteq \{1\l n\}$ with $|I|=|J|=k$ such that $\bigcup_{i\in I}A_i$
is $G$-complete with $\bigcup_{j\in J} B_j$.
\end{itemize}
\end{thm}

A {\em $k$-clique} means a clique of cardinality $k$. 
If $X$ is a clique of a digraph $G$, a vertex in $X$ is a {\em source} of $X$ if it is $G$-adjacent to every other vertex in $X$, and a 
{\em sink} if it is $G$-adjacent from every other vertex in $X$.
If $k,m\ge 1$ are integers, a vertex $v$ of a digraph $G$ is {\em $(k,m)$-rich} if there exist $k$ pairwise disjoint $m$-cliques
$A_1\l A_k\subseteq N^+(v)$, and $k$ pairwise disjoint $m$-cliques $B_1\l B_k\subseteq N^-(v)$, such that 
$A_1\cup\cdots\cup A_k$ is $G^*$-complete with $B_1\cup \cdots\cup B_k$.

\begin{thm}\label{outnbrs}
For all integers $k,m\ge 1$ there exists $t$ with the following property. Let $G$ be a $\lambda$-spread digraph 
such that
no vertex of $G$ is $(k,m)$-rich.
Then $V(G)$ can be partitioned into $t$ sets $X_1\l X_t$ such that for $1\le i\le t$,
either no $(m+1)$-clique of $G[X_i]$ has a source or no $(m+1)$-clique of $G[X_i]$ has a sink.
\end{thm}
\Proof
Choose $n$ such that \ref{bigramsey} holds, and let $t=4nm$. We claim that $t$ satisfies \ref{outnbrs}.
Let $G$ be as in the theorem. Let $P$ be the set of vertices of $G$ such that there do not exist $n$ 
pairwise disjoint $m$-cliques in $N^+(v)$, and let $Q$ be the set such that there do not exist $n$ 
pairwise disjoint $m$-cliques in $N^-(v)$. Suppose first that some vertex $v$ belongs to neither of $P,Q$. Then 
there exist $n$ pairwise disjoint $m$-cliques
$A_1\l A_n\subseteq N^+(v)$, and there exist $n$ pairwise disjoint $m$-cliques $B_1\l B_n\subseteq N^-(v)$.
Since $G$ is $\lambda$-spread, 
\ref{bigramsey} implies that there exist $I, J\subseteq \{1\l n\}$ with $|I|=|J|=k$ 
such that $\bigcup_{i\in I}A_i$
is $G^*$-complete with $\bigcup_{j\in J} B_j$, that is, $v$ is $(k,m)$-rich, a contradiction. This proves
that $P\cup Q=V(G)$.

For each vertex $v\in P$, choose a maximal set of pairwise disjoint $m$-cliques included in $N^+(v)$, and let
the union of the members of this set be $P_v$. Then each $|P_v|<
nm$, and has nonempty intersection with every $m$-clique included in $N^+(v)$. Let $H$ be the digraph
with vertex set $P$ and edge set consisting of the edges with tail $v$ and head in $P_v$, for each $v\in P$. Then every vertex of $H$
has outdegree less than $nm$, and so $\chi(H)\le 2nm$. Let $X$ be a stable set of $H^*$. It follows
that for each $v\in X$, there is no $m$-clique included in $N^+(v)\cap X$ (since $P_v$ has nonempty intersection with 
every such clique, and $P_v\cap X=\emptyset$ because $X$
is stable in $H^*$). Consequently there is no $(m+1)$-clique
included in $X$ that has a source. But $P$ can be partitioned into $\chi(H)\le 2nm=t/2$ such sets $X$, and similarly
we can partition $Q$. This proves \ref{outnbrs}.~\bbox

Choose a function $\phi$ such that for all $k,m\ge 0$, 
setting $t=\phi(k,m)$ satisfies \ref{outnbrs}. Let $\phi$ be fixed for the remainder of this section.

Next we prove \ref{stars} for acyclic digraphs (a digraph is {\em acyclic} if it has no directed cycle).
\begin{thm}\label{acyclicstars}
There exists $c_0$ such that $\chi(G)\le c_0$ for every acyclic $\lambda$-spread digraph $G$
with $\omega(G)\le \kappa$.
\end{thm}
\Proof
Let $t=\phi(1,\kappa-1)$, and let
$c_0= t\tau$. We claim that $c_0$ satisfies the theorem.

Let $G$ be a $\lambda$-spread acyclic digraph with $\omega(G)\le \kappa$.
Now no vertex of $G$ is $(1,\kappa-1)$-rich, because then $G$ would have a clique of cardinality $2\kappa-1>\kappa$.
By \ref{outnbrs}, $V(G)$ can be partitioned into $t$ sets $X_1\l X_t$, such that for $1\le i\le t$, either 
no $\kappa$-clique of $G[X_i]$ has a source or no $\kappa$-clique of $G[X_i]$ has a sink. But every $\kappa$-clique
has both a source and a sink, since $G$ is acyclic, and so $\omega(G[X_i])<\kappa$, and consequently $\chi(X_i)\le \tau$. Hence
$\chi(G)\le t\tau=c_0$. This proves \ref{acyclicstars}.~\bbox

A digraph $G$ is {\em $(h,k)$-out-orderable} if there is a partition $X_1\l X_n$ of its vertex set, such that
for $1\le i\le n$, $\chi(X_i)\le h$, and each vertex of $X_i$ has at most $k-1$ out-neighbours in $X_{i+1}\cup \cdots\cup X_n$.
We define $(h,k)$-in-orderable similarly. If $G$
is a digraph, we say $X\subseteq V(G)$ is {\em acyclic} if $G[X]$ is acyclic.

\begin{thm}\label{outorderable}
If the digraph $G$ is $(h,k)$-out-orderable, then there is a partition of $V(G)$ into $hk$ acyclic sets.
\end{thm}
\Proof
Let $X_1\l X_n$ be as in the definition of $(h,k)$-out-orderable.
Let $J$ be the graph with vertex set $V(G)$ in which $u,v$ are $J$-adjacent if $u$ is $G$-adjacent to $v$ and $i\le j$
where $u\in X_i$ and $v\in X_j$.
Since $J[X_i]$ is $h$-colourable for each $i$, there is a partition $Y_1\l Y_h$ of $V(J)$ such that $X_i\cap Y_j$ is stable
for $1\le i\le n$ and $1\le j\le h$. But every nonempty induced subgraph of $J[Y_j]$ has a vertex with degree in $J$ less than $k$
(choose a vertex of $Y_j$ in $X_i$ for the smallest $i$ with $X_i\cap Y_j$ nonempty); and so $J[Y_j]$ is $(k-1)$-degenerate and
hence $k$-colourable. Consequently $\chi(J)\le hk$. But for each stable set $Y$ of $J$, $G[I]$ is acyclic. 
This proves \ref{outorderable}.~\bbox

A digraph $G$ is {\em $(h,k)$-robust} if for every nonempty subset $Z\subseteq V(G)$ with $\chi(Z)\le h$, some vertex of $Z$
has at least $k$ out-neighbours in $V(G)\setminus Z$ and at least $k$ in-neighbours in $V(G)\setminus Z$. 

\begin{thm}\label{robustpartition}
Let $h,k\ge 0$; then for every digraph $G$ there is a partition of $V(G)$ into three sets $P,Q,R$
such that $G[P]$ is $(h,k)$-out-orderable, $G[Q]$ is $(h,k)$-in-orderable and $G[R]$ is $(h,k)$-robust.
\end{thm}
\Proof
We proceed by induction on $|V(G)|$. If $G$ is $(h,k)$-robust we are done, so we may assume that there is a nonempty subset
$Z\subseteq V(G)$ with $\chi(Z)\le h$, such that for each $v\in Z$, either $|N^+(v)\setminus Z|<k$ or
$|N^-(v)\setminus Z|<k$. Let $X_1$ be the set of vertices $v\in Z$ such that $|N^+(v)\setminus X|<k$,
and $Y_1=Z\setminus X_1$. From the inductive hypothesis there is a partition $P,Q,R$ of $V(G)\setminus Z$
such that $G[P]$ is $(h,k)$-out-orderable, $G[Q]$ is $(h,k)$-in-orderable and $G[R]$ is $(h,k)$-robust. Let 
$X_2\l X_n$ be a partition of $P$ such that 
for $2\le i\le n$, $\chi(X_i)\le h$, and each vertex of $X_i$ has at most $k-1$ out-neighbours in $X_{i+1}\cup \cdots\cup X_n$.
Then the sequence $X_1\l X_n$ shows that $G[P\cup X_1]$ is $(h,k)$-out-orderable. Similarly $G[Q\cup Y_1]$ is
$(h,k)$-in-orderable, and so the partition $P\cup X_1,Q\cup Y_1,R$ satisfies the theorem. This proves \ref{robustpartition}.~\bbox

We recall that $\kappa, \lambda$ and $\tau$ are fixed integers satisfying $\kappa\ge 2$, $\lambda\ge 0$, 
and $\chi(J)\le \tau$ for
every $\lambda$-spread digraph $J$ with $\omega(J)<\kappa$.
Let us define $\Lambda = 2\lambda^2+\lambda$ (throughout the remainder of this section).
\begin{thm}\label{userobust}
Let $G$ be a $\lambda$-spread digraph and let $X\subseteq V(G)$ be nonempty.
If every vertex in $X$ has at least $\Lambda$ 
out-neighbours in $X$
and at least $\Lambda$ in-neighbours in $X$, then $G$ is not $(|X|\tau,|X|+\Lambda)$-robust.
\end{thm}
\Proof
We claim first:
\\
\\
(1) {\em For each vertex $v$ of $G$, if $A\subseteq N^+(v)$ and $B\subseteq N^-(v)$ with $|A|=\lambda$, then some vertex of $A$
is $G^*$-adjacent with at least $|B|/\lambda-1$ members of $B$.}
\\
\\
There are fewer than $\lambda$ members of $B$ that have no $G^*$-neighbour in $A$, since 
$G$ is $\lambda$-spread. So
all the others have at least one $G^*$-neighbour in $A$; and so some vertex in $A$ is $G^*$-adjacent 
with at least $(|B|-\lambda)/\lambda$ of them.
This proves (1).

\bigskip
Now let $X\subseteq V(G)$ be nonempty, such that every vertex in $X$ has at least $\Lambda$
out-neighbours in $X$
and at least $\Lambda$ in-neighbours in $X$.
Let $P$ be the set of vertices not in $X$ with at least $2\lambda$ $G^*$-neighbours in $X$.
\\
\\
(2) {\em For each $u\in X$, $u$ is $G^*$-adjacent with fewer than $2\lambda$
vertices in $V(G)\setminus (P\cup X)$.}
\\
\\
For suppose not; then from the symmetry we may assume that there is a set $A$ of 
in-neighbours of $u$ in $V(G)\setminus (P\cup X)$ with $|A|=\lambda$. But $u$ has at least $\Lambda$ 
out-neighbours in $X$; and so by (1),
some vertex in
$V(G)\setminus (P\cup X)$ has at least $2\lambda$ neighbours in $X$, and therefore belongs to $P$, a contradiction.
This proves (2).

\bigskip

Suppose that
there exists $v\in P$ with at least $|X|+\Lambda$ 
out-neighbours in $V(G)\setminus P$
and at least $|X|+\Lambda$ in-neighbours in $V(G)\setminus P$. Since $v$ has at least $2\lambda$ $G^*$-neighbours
in $X$, from the symmetry we may assume that $v$ has at least $\lambda$ out-neighbours in $X$. Let $Y$ be the
set of vertices in $V(G)\setminus (P\cup X)$ that are in-neighbours of $v$. Then $|Y|\ge \Lambda$. 
Since $v$ has at least $\lambda$
out-neighbours in $X$, (1) implies that one of these out-neighbours, say $u$, is $G^*$-adjacent with at least 
$|Y|/\lambda-1\ge 2\lambda$
vertices in $Y$, contrary to (2). Thus there is no such $v$. But 
$\chi(P)\le |X|\tau$ since every vertex in $P$ has a neighbour in $X$. If $P\ne \emptyset$ it follows 
that $G$ is not $(|X|\tau,|X|+\Lambda)$-robust, as required, so we may assume that $P=\emptyset$. But then by (2),
every vertex in $X$ is $G^*$-adjacent with fewer than $2\lambda$ vertices in $V(G)\setminus X$, and so again, 
$G$ is not $(|X|\tau,|X|+\Lambda)$-robust.
This proves \ref{userobust}.~\bbox

If $G$ is a digraph and $u,v,w$ are vertices, pairwise $G^*$-adjacent, such that one of them is $G$-adjacent from the other two,
we call $\{u,v,w\}$ a {\em transitive triangle}.
Next we need:

\begin{thm}\label{gettri}
There exists $k_0$ with the following property.
Every non-null $(3\Lambda\tau,k_0)$-robust $\lambda$-spread digraph has a transitive triangle.
\end{thm}
\Proof
By the bipartite version of Ramsey's theorem~\cite{beineke},
for all $n\ge 0$ there exists $f(n)\ge n$ such that 
for every partition of the edges of the complete bipartite graph $K_{f(n),f(n)}$ into two sets, 
either the first set includes the edge set of a $K_{n,n}$ 
subgraph, or the second set includes the edge set of a $K_{\lambda,\lambda}$-subgraph.
Let $k_0=f(f(f(\Lambda)))$.
Suppose that $G$ is a $(3\Lambda\tau,k_0)$-robust $\lambda$-spread digraph with no transitive triangle.
Since $G$ is $(3\Lambda\tau,k_0)$-robust, every vertex has at least $k_0$ out-neighbours and $k_0$ in-neighbours. 
Let $v\in V(G)$; then since $G$ is $\lambda$-spread, there exist $A_1\subseteq N^+(v)$ and $B_1\subseteq N^-(v)$
with $|A_1|=|B_1|=f(f(\Lambda))$ such that $A_1$ is $G^*$-complete with $B_1$; and since there is no transitive triangle it follows
that every vertex in $A_1$ is $G$-adjacent to every vertex in $B_1$. 
Choose $a\in A_1$. Since $a$ has at least $k_0$ in-neighbours, and none of them belong to $A_1\cup B_1$ (because $A_1$ is stable since
there is no transitive triangle)
there is a set $C$ of vertices in $V(G)\setminus (A_1\cup B_1)$ all $G$-adjacent to $a$, with $|C|=k_0$.
Since $a$ is $G$-adjacent to every vertex in $B_1$, and $|B_1|,|C|\ge f(f(\Lambda))$, and there is no transitive triangle, 
there exist $B_2\subseteq B_1$ and $C_1\subseteq C$ with $|B_2|=\Lambda$ and $|C_1|=f(\Lambda)$ such that every vertex in $B_2$ is $G$-adjacent 
to every vertex
in $C_1$. Choose $b\in B_2$. Since $b$ is $G$-adjacent to every vertex in $C_1$ and from every vertex in $A_1$, and $|A_1|,|C_1|\ge f(\Lambda)$,
there exist $A_2\subseteq A_1$ and $C_2\subseteq C_1$ with $|A_2|,|C_2|=\Lambda$ such that $A_2$ is $G$-complete from
$C_2$. Since $|A_2|, |B_2|,|C_2|= \Lambda$, every vertex in $A_2\cup B_2\cup C_2$ has at least $\Lambda$ out-neighbours
and $\Lambda$ in-neighbours in $A_2\cup B_2\cup C_2$, contrary to \ref{userobust} (taking $X=A_2\cup B_2\cup C_2$). This proves~\ref{gettri}.~\bbox

A tournament $H$ is {\em regular} if all its vertices have the same outdegree, and they all have the same indegree;
and it follows that $|V(H)|$ is odd, $|V(H)|=2m+1$ say, and all vertices have indegree and outdegree $m$. 
A tournament $H$ is  {\em cyclic} if it has an odd number of vertices, say $2m+1$, and its vertex set can be ordered as 
$\{v_1\l v_{2m+1}\}$ such that for $1\le i<j\le 2m+1$, $v_i$ is $H$-adjacent to $v_j$ if and only if $j-i\le m$.
\begin{thm}\label{cyclic}
Let $H$ be a regular tournament with $2m+1$ vertices, and let $v\in V(H)$; and suppose there is no directed cycle
with vertices $p\d q\d r\d s\d p$ in order such that $p,r$ are out-neighbours of $v$ and $q,s$ are in-neighbours of $v$.
Then $H$ is cyclic.
\end{thm}
\Proof
Let $J$ be the subdigraph of $H$ with vertex set $V(H)$ and edge set consisting of all edges between $N^+(v)$ and $N^-(v)$. If $J$ has a 
directed cycle, take the shortest such directed cycle $C$; then $C$ is induced and so has length four, a contradiction. Thus $J$ has no 
directed cycle, and so $V(H)\setminus \{v\}$ can be ordered as $\{v_1\l v_{2m}\}$ such that for every edge of $J$, its tail is earlier
than its head. Thus $v_1$ has $m$ out-neighbours in $J$, and since it has only $m$ out-neighbours in $H$, it follows that 
$v_1\in N^+(v)$, and $v_1$ is adjacent from every other vertex in $N^+(v)$. If $v_2\in N^+(v)$, then it has $m$ out-neighbours
in $N^-(v)$, and $v_1$ is another, a contradiction; so $v_2\in N^-(v)$, and hence $v_2$ is adjacent to every vertex in $N^+(v)$
except $v_1$, and therefore adjacent from every other vertex in $N^-(v)$. More generally,
we claim that 
\begin{itemize}
\item for $i$ odd, $v_i\in N^+(v)$ and $v_i$ is $H$-adjacent from all vertices of $N^+(v)$ except $v_1,v_3\l v_{i-2}$; and
\item for $i$ even, $v_i\in N^-(v)$ and $v_i$ is $H$-adjacent from all vertices of $N^-(v)$ except $v_2,v_4\l v_{i-2}$.
\end{itemize}
We prove this claim by induction on $i$. Suppose then that it holds for all smaller values of $i$, and first suppose 
that $v_i\in N^+(v)$. Then $v_i$ is $H$-adjacent to $v_h$ for all odd $h<i$ (from the inductive hypothesis applied to $v_h$),
and there are $\lfloor i/2\rfloor$ such values of $h$. 
Also, there are exactly $\lceil i/2\rceil -1$ values of $j$ with $j<i$ such that $v_j\in N^-(v)$, from the inductive hypothesis,
and $v_i$ is $H$-adjacent to each of the remaining $m+1-\lceil i/2\rceil$ vertices in $N^-(v)$, from the property of the ordering.
Since 
$v_i$ has outdegree exactly $m$ in $H$, it must be the case that $\lfloor i/2\rfloor + m+1-\lceil i/2\rceil\le m$. So $i$ is odd,
and moreover $v_i$ has no further outneighbours; and so $v_i$ is $H$-adjacent from
all vertices of $N^+(v)$ except $v_1,v_3\l v_{i-2}$ as claimed. 

Now suppose that $v_i\in N^-(v)$.
Thus $v_i$ is $H$-adjacent to $v_h$ for all even $h<i$ by the inductive hypothesis, and to every vertex of $N^+(v)$
except for the vertices $v_j$ with $j$ odd and $j<i$, because of the ordering. 
There are $\lceil i/2\rceil-1$ outneighbours of the first kind, and 
$m-\lfloor i/2\rfloor$ of the second kind; and in addition $v_i$ is $H$-adjacent to $v$. Consequently
$\lceil i/2\rceil-1 + m-\lfloor i/2\rfloor +1\le m$, and so $i$ is even, and 
$v_i$ is $H$-adjacent from every vertex of $N^-(v)$
except $v_2,v_4\l v_{i-2}$. This proves the inductive statement, and so proves
\ref{cyclic}.~\bbox

\begin{thm}\label{sinks}
There exist $k_1,c_1\ge 0$ with the following property. 
Let $G$ be a $(4\Lambda\tau,k_1)$-robust $\lambda$-spread digraph with 
$\omega(G)\le \kappa$, such that no $(\lceil \kappa/2\rceil+1)$-clique of $G$ has a source. Then $\chi(G)\le c_1$.
\end{thm}
\Proof
Let $k_0$ satisfy \ref{gettri}, and let $k_1=\max(k_0, 5\Lambda)$.
Choose $n\ge 0$ such that
for every partition of the edges of the complete bipartite graph $K_{n,n}$ into two sets,
either the first set includes the edge set of a $K_{\Lambda,\Lambda}$
subgraph, or the second set includes the edges of a $K_{\lambda,\lambda}$-subgraph.

Let $m=\lfloor \kappa/2\rfloor$. Choose $k$ such that for every partition of the edges
of $K_{k,k}$ into $m^4+1$ sets, one of the sets includes all the edges of some $K_{n,n}$ subgraph.
Let $c_1=\phi(k,m)\tau$. 

Now let $G$ be as in the theorem. We claim that $\chi(G)\le c_1$. 
If $\kappa$ is even then since no clique of $G$ has a vertex of outdegree $\kappa/2$ (because 
no $(\lceil \kappa/2\rceil+1)$-clique of $G$ has a source), it follows that $\omega(G)<\kappa$ and so $\chi(G)\le \tau\le c_1$.
We may therefore assume that $\kappa$ is odd, and so $\kappa = 2m+1$.
If $\kappa=3$, then since no $(\lceil \kappa/2\rceil+1)$-clique of $G$ has a source, it follows that $G$ has no transitive triangle,
contrary to \ref{gettri}. Thus $\kappa>3$ and so $m\ge 2$.
\\
\\
(1) {\em No vertex of $G$ is $(k,m)$-rich.}
\\
\\
Because suppose that 
$v$ say is $(k,m)$-rich. Let $A_1\l A_k, B_1\l B_k$
be cliques as in the definition of $(k,m)$-rich. For $1\le i\le k$ let $A_i=\{a_1^i\l a_m^i\}$
and $B_i=\{b_1^i\l b_m^i\}$, choosing the numbering such that if $A_i$ is a transitive tournament, then $a_p^i$ is
$G$-adjacent to $a_q^i$ for all $p<q$, and if $B_i$ is transitve then $b_p^i$ is $G$-adjacent to $b_q^i$ for all $p<q$.
For $1\le i,j\le m$, if there is a directed cycle of length four, with vertices
$a_p^i\d b_q^j\d a_r^i\d b_s^j\d a_p^i$ in order, we say the pair $(i,j)$ has {\em type $(p,q,r,s)$} (choosing some such quadruple
arbitrarily if there is more than one), and type 0 otherwise. By the choice of $k$, we may assume that
all the pairs $(i,j)$ for $1\le i,j\le n$ have the same type. If this type is nonzero, say $(p,q,r,s)$, let $X$ be the set
$$\bigcup_{1\le i\le \Lambda}\{a_p^i,a_r^i\}\cup \bigcup_{1\le j\le \Lambda}\{b_q^j,b_s^j\}.$$
Every vertex in $X$ has at least $\Lambda$
out-neighbours and $\Lambda$ in-neighbours in $X$, and $|X|=4\Lambda$, contrary to \ref{userobust}.

Thus for $1\le i,j\le n$, $(i,j)$ has type 0. From \ref{cyclic}, $G[A_i\cup B_j\cup\{v\}]$ is cyclic, 
and therefore both $A_i,B_j$ are transitive tournaments,
and from the choice of numbering, $a_p^i$ is $G$-adjacent to $a_q^i$ for all $p<q$, and similarly
$b_p^j$ is $G$-adjacent to $b_q^j$ for all $p<q$. Since $G[A_i\cup B_j\cup\{v\}]$
is cyclic, it follows that for $1\le p,q\le m$, $a_p^i$ is $G$-adjacent to $b_q^j$ if and only if $q\le p$.

Suppose that for some $p,q\in \{1\l n\}$, $a_1^p$ is $G$-adjacent from $a_m^q$. Then the subdigraph induced on
$\{a_1^p,a_m^q\}\cup B_1$ is an $(m+2)$-clique with a source (namely $a_m^q$), a contradiction.
Now $b_1^1$ is $G$-adjacent to each of $a_1^1\l a_1^n$ and $G$-adjacent from each of $a_m^1\l a_m^n$. Consequently
since $G$ is $\lambda$-spread, from the definition of $t$ there exist $A\subseteq \{a_1^1\l a_1^n\}$ and $C\subseteq \{a_m^1\l a_m^n\}$
with $|A|=|C|=\Lambda$, such that $A$ is $G$-complete to $C$. Let $B=\{b_1^1\l b_{\Lambda}^1\}$; then
$B$ is $G$-complete to $A$, and $C$ is $G$-complete to $B$.
Consequently, every vertex in $A\cup B\cup C$ has at least $\Lambda$ out-neighbours and $\Lambda$ in-neighbours in this set.
But this contradicts \ref{userobust}. This proves (1).

\bigskip

By \ref{outnbrs}, $V(G)$ can be partitioned into $\phi(k,m)$ subsets  such that
for each such subset $Y$ say, either no $(m+1)$-clique of $G[Y]$ has a source, or none has a sink. In either case it follows
that $\omega(G[Y])<\kappa$ and so $\chi(Y)\le \tau$ and hence $\chi(G)\le \phi(k,m)\tau=c_1$. This proves \ref{sinks}.~\bbox

\noindent{\bf Proof of \ref{stars}.\ \ }
As discussed at the beginning of this section, it suffices to show that for some $c\ge 0$, $\chi(G)\le c$ for 
every $\lambda$-spread digraph $G$ 
with $\omega(G)\le \kappa$. 
Let $c_0$ satisfy \ref{acyclicstars}, let $c_1,k_1$ satisfy \ref{sinks}, 
and let $c=4\Lambda\tau k_1 c_0 + \phi(1,\lceil \kappa/2\rceil)c_1$.

Now let $G$ be a $\lambda$-spread digraph $G$
with $\omega(G)\le \kappa$.
By \ref{robustpartition} and \ref{outorderable} there is a subset $R\subseteq V(G)$
such that $G[R]$ is $(4\Lambda\tau,k_1)$-robust and $V(G)\setminus R$ can be partitioned into $4\Lambda \tau k_1$ acyclic sets; and each
of the latter induces a $c_0$-colourable digraph by \ref{acyclicstars}. Thus $\chi(G)\le 4\Lambda \tau k_1 c_0 +\chi(R)$, 
so it remains to bound $\chi(R)$.

No vertex is $(1,\lceil \kappa/2\rceil)$-rich, since that would imply that $G$ contains a $(\kappa+1)$-clique. By \ref{outnbrs},
$V(G)$ can be partitioned into $\phi(1,\lceil \kappa/2\rceil)$ subsets such that for each such subset $Y$ say, either no
$(\lceil \kappa/2\rceil+1)$-clique of $G[Y]$ has a source or none has a sink.
From \ref{sinks} it follows that $\chi(Y\cap R)\le c_1$ for each $Y$, and so $\chi(R)\le \phi(1,\lceil \kappa/2\rceil)c_1$. 
Consequently $\chi(G)\le 4\Lambda \tau k_1 c_0 + \phi(1,\lceil \kappa/2\rceil) c_1=c $.
This proves \ref{stars}.~\bbox

\end{document}